\documentclass[11pt]{amsart}

\usepackage{geometry}
\usepackage{graphicx}
\usepackage{bm,physics}

\renewcommand{\qq}{\bm{q}}

\newcommand{\vv}{\bm{v}}
\newcommand{\nn}{\bm{n}}
\newcommand{\E}{\mathcal{E}}
\newcommand{\T}{\mathcal{T}}
\newcommand{\dT}{\partial\T}
\newcommand{\intelem}[2][h]{\left(#2\right)_{\T_{#1}}}
\newcommand{\intfa}[2][h]{\langle#2\rangle_{\dT_{#1}}}

\theoremstyle{definition}

\title{A non-overlapping Schwarz algorithm for the HDG method}
\author{Issei OIKAWA}

\begin{document} 
\begin{abstract}
In this paper, we present two non-overlapping Schwarz algorithms for the hybridizable discontinuous Galerkin (HDG) method. 
The first algorithm is based on the Neumann-Neumann method. 
The second one is an iterative algorithm uses both trace and flux interface unknowns on interfaces 
between subdomains.
Numerical results are provided to verify the validity of our algorithms.  
\end{abstract}

\maketitle

\section{Introduction}
Let $\Omega \subset \mathbb{R}^d~(d = 2,3)$ be a bounded polygonal or polyhedral domain. 
We consider the Poisson equation as a model problem:
\begin{subequations}\label{pois}
\begin{align}
-\Delta u &= f \qquad \text{ in } \Omega,  \label{pois1}\\ 
u &= 0 \qquad  \text{ on }\partial\Omega, \label{poisd}
\end{align}
\end{subequations}
where $f$ is a given function.
Let $\Omega = \Omega_1 \cup \Omega_2$ and $\Omega_1 \cap \Omega_2 = \emptyset$. 
Let $\Gamma_{12}$ denote the interface between the subdomains.
The problem \eqref{pois} can be rewritten into 
\begin{subequations}\label{subpois}
\begin{align}
-\Delta u_i &= f \text{ in } \Omega_i \quad (i=1,2),  \\ 
u_i &= 0 \text{ on } \partial\Omega_i \cap \partial\Omega \quad (i=1,2),  \\
\pdv{u_1}{n_1} + \pdv{u_2}{n_2} &= 0 \text{ on } \Gamma_{12},  \\ 
u_1 &= u_2 \text{ on } \Gamma_{12}.
\end{align}
\end{subequations}
Here $n_i$ is the outer unit normal vector to $\partial\Omega_i$. 
A non-overlapping Schwarz algorithm (cf. \cite{ToWi2005,VJN2015}) is a family of 
domain decomposition methods to solve the subproblems separately and is widely used to compute numerical solutions in parallel. 
The optimized Schwarz method is introduced by Lions \cite{LionsOSM1}, which is based on Robin interface condition and can be applied both overlapping and non-overlapping cases. 
In \cite{GaHa2014,GaHa2015,GaHa2018}, the optimized Schwarz method  of the HDG method is proposed and analyzed.
The Neumann-Neumann method \cite{Mandel1993} and  the FETI (or Dirichlet-Dirichlet) method \cite{FaRo1991} are also well known as a non-overlapping algorithm, however, there is no application to the HDG method to the best of our knowledge.

In this paper, we present two non-overlapping Schwarz algorithms intended to apply to the hybridizable discontinuous Galerkin (HDG) method \cite{CGL2009}.
The first algorithm is a direct application of the Neumann-Neumann algorithm to the HDG method.
In the HDG method, the problem is split into many elements, and numerical trace and flux on inter-element boundaries are introduced. 
The Neumann-Neumann and the FETI methods introduce an interface unknown on interfaces between subdomains. 
The interface unknown plays a similar role to the numerical trace or flux that reduces the jump of the solution on interface, 
so the Neumann-Neumann and the FETI methods are highly compatible with the HDG method.
We also another non-overlapping algorithm. 
The key idea is to use alternatively trace- and flux-interface unknowns on interfaces between subdomains. 
The interface unknowns are updated by using the numerical trace or flux of the solutions on subdomains. 
Therefore, the second algorithm is different from the other algorihtms. We note that the proposed algorithm 
does not include any parameter in iteration.

The rest of the paper is organized as follows.
In Section 2, we describe the Neumann-Neumann algorithm and application to the HDG method, i.e. the first algorithm, in two-subdomain case.
In Section 3, the second algorithm and discretization by the HDG method are presented. 
Numerical results are also provided to verify the validity of the proposed algorithm.
All numerical computation were carried out by FreeFEM \cite{MR3043640} and Julia \cite{julia}.

\section{The Neumann-Neumann algorithm}
\subsection{The Neumann-Neumann algorithm}
We begin by recalling the Neumann-Neumann algorithm.
\begin{enumerate}
\item[Step 1.] Set the initial value: $u_{\Gamma}^0$
\item[Step 2.] Repeart Steps 3--5 for $n\ge 0$ until convergence
\item[Step 3.] Compute $u_i^{n+1/2}$ $(i=1,2)$ by solving
\begin{align*}
-\Delta u_i^{n+1/2} &= f \text{ in } \Omega_i,  \\ 
u_i^{n+1/2} &= 0 \text{ on } \partial\Omega_i \cap \partial\Omega,  \\
u_i^{n+1/2} &= u_\Gamma^n \text{ on } \Gamma_{12}. 
\end{align*}
\item[Step 4.] 
Compute $\psi_i^{n+1} ~ (i=1,2)$ by solving
\begin{align*}
-\Delta \psi_i^{n+1} &= 0  \text{ in } \Omega_i, \\
\psi_i^{n+1} &= 0 \text{ on } \partial\Omega_i \cap \partial\Omega, \\
\pdv{\psi_i}{n_i}^{n+1} &= 
\pdv{u_1}{n_1}^{n+1/2}
+  \pdv{u_2}{n_2}^{n+1/2}  \text{ on } \Gamma_{12}.
\end{align*}
\item[Step 5.] Update 
\[
u_\Gamma^{n+1} =  u_\Gamma^{n} -\theta(\psi_1^{n+1}+
\psi_2^{n+1}) \text{ on } \Gamma_{12}. 
\]
Here $\theta>0$ is a constant parameter.
\end{enumerate}

\subsection{The HDG approximation}
Let $\T_{h}$ be a mesh of $\Omega$. The set of all edges of $K \in \T_h$ is denoted by $\E_h$.
Define $\T_{ih} = \{ K \in \T_h : K \subset \Omega_i \}$ 
and $\E_{ih} = \{ e \in \E_h : e \subset \Omega_i \}$ for $i = 1,2$.
Let $\mathcal{E}_{ih}$ denote the set of all edges of $\partial K$ for $K \in \mathcal{T}_{ih}$. 
We assume that $\Gamma_{12} = \bigcup_{e \subset \Gamma_{12}, e \in \E_h} e$.
We introduce finite dimensional spaces $\bm{V}(K)$, $W(K)$ and $M(e)$ 
to approximate $\qq|_K, u|_K$ and $u|_e$, respectively, where $K \in \mathcal{T}_h$  and $e \in \mathcal{E}_h$.
The global approximate spaces are constructed as 
\begin{align*}
\bm{V}_{ih} &= \{ \vv \in L^2(\Omega_i)^d : \vv|_K \in \bm{V}(K) ~ \forall K \in \mathcal{T}_{ih}\}, \\ 
W_{ih} &= \{ w \in L^2(\Omega_i) : w|_K \in W(K) ~ \forall K \in \mathcal{T}_{ih}\}, \\ 
M_{ih} &= \{ \mu \in L^2(\mathcal{E}_{ih}) : \mu|_e \in M(e) ~ \forall e \in \mathcal{E}_{ih}, \quad \mu|_{\partial\Omega} = 0\}. 
\end{align*}
The inner product on a domain $D$ or a curve $F$ is denoted as 
\[
(u,w)_D = \int_D uw dx, \quad \intfa[F]{\lambda, \mu} = \int_F \lambda\mu ds.    
\] 
The piecewise inner products are denoted as
\begin{align*}
&\intelem[ih]{\qq,\vv} = \sum_{K \in \T_{ih}} \int_K \qq\cdot\vv dx, \quad 
\intelem[ih]{u,w} = \sum_{K \in \T_{ih}} \int_K uw dx, \\ 
&\intfa[ih]{\lambda,\mu} = \sum_{K \in \T_{ih}}
\int_{\partial K} \lambda\mu ds. 
\end{align*}
The HDG approximation of the Neumann-Neumann algorithm is as follows.
\begin{enumerate}
\item[Step 1.] Set the initial value $\widehat{u}_{\Gamma}^0$.
\item[Step 2.] Repeat Steps 3-5 for $n \ge 0$ until convergence.
\item[Step 3.]Solve the following equations  to get $(\qq_i^{n+1/2},u_i^{n+1/2}, \widehat{u}_i^{n+1/2}) \in 
\bm{V}_{ih} \times W_{ih} \times M_{ih}$ for $i=1,2$: 
\begin{align*}
\intelem[ih]{\qq_i^{n+1/2}, \vv} - \intfa[ih]{u_i^{n+1/2}, \nabla\cdot \vv} 
+ \intfa[ih]{\widehat{u}_h^{n+1/2}, \vv\cdot\nn} & = 0 && \forall \vv \in \bm{V}_{ih} \\
-\intelem[ih]{\qq_i^{n+1/2}, \nabla w} + \intfa[ih]{\widehat{\qq}_i^{n+1/2}\cdot\nn, w}  
&= (f,w)_{\Omega_i} && \forall w \in W_{ih}, \\ 
\intfa[ih]{\widehat{\qq}_i^{n+1/2}\cdot\nn, \mu} & = 0 && \forall \mu \in M_{ih}, 
\end{align*}
where   
\begin{align*}
\widehat{\qq}_i^{n+1/2} \cdot\nn &= \qq_i^{n+1/2} \cdot \nn + \tau(u_i^{n+1/2} - \widehat{u}_i^{n+1/2}), \\ 
\widehat{u}_i^{n+1/2} &= \widehat{u}_\Gamma^n \quad \text{ on } \Gamma_{12}.
\end{align*}

\item[Step 4.]
Solve the following equations to get    
$(\bm\sigma_i^{n+1},\xi_i^{n+1}, \widehat{\xi}_i^{n+1}) \in 
\bm{V}_{ih} \times W_{ih} \times M_{ih}$ for $i=1,2$:
\begin{align*}
\intelem[ih]{\bm\sigma_{i}^{n+1}, \vv} - \intfa[ih]{\xi_{i}^{n+1}, \nabla\cdot \vv} 
+ \intfa[ih]{\widehat{\xi}_{i}^{n+1}, \vv\cdot\nn} & = 0 && \forall \vv \in V_{ih}, \\
-\intelem[ih]{\widehat{\bm{\sigma}}_{i}^{n+1}, \nabla w} + \intfa[ih]{\widehat{\bm{\sigma}}_{i}^{n+1}\cdot\nn, w}  
&= 0 && \forall w \in W_{ih}, \\ 
\intfa[ih]{\widehat{\bm{\sigma}}_{i}^{n+1}\cdot\nn, \mu}  = 
\intfa[\Gamma_{12}]{\widehat{\qq}_{1}^{n+1/2}\cdot\nn_1 + \widehat{\qq}_{2}^{n+1/2}\cdot\nn_2,\mu}& && \forall \mu \in M_{ih}, 
\end{align*}
where
\[
\widehat{\bm\sigma}_{i}^{n+1} \cdot\nn = \bm\sigma_{i}^{n+1} \cdot \nn + \tau(\xi_{i}^{n+1} - \widehat{\xi}_{i}^{n+1})
\]
and $\tau>0$ is a stabilization parameter.
\item[Step 5.] 
Update 
\[
\widehat{u}^{n+1}_\Gamma = \widehat{u}_{\Gamma}^{n} + 
\theta(\widehat{\xi}_{1h}^{n+1} + \widehat{\xi}_{2h}^{n+1}).
\]
\end{enumerate}

\subsection{Numerical results} 

We consider the following test problem:
\begin{subequations} \label{testprob}
\begin{align} 
-\Delta u &= 2\pi^2\sin(\pi x)\sin(\pi y) 
&&  \text{ in } \Omega := (0.1)^2, \\ 
u &= 0 && \text{ on } \partial\Omega.
\end{align}
\end{subequations}
The domain is decomposed into 
$\Omega_1 = (0,1/2)\times(0,1)$ and $\Omega_2 = (1/2,1)\times(0,1)$.
We use unstructured meshes for each subdomain, where there is no hanging node on the interface.   
All approximation spaces are piecewise or edgewise polynomials of degree one. 
The stabilization parameter of the HDG method is given by $\tau \equiv 1$.

\begin{figure}[h]
\includegraphics[width=300pt]{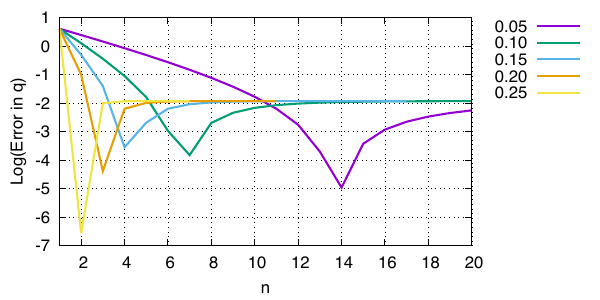}  \\
\includegraphics[width=300pt]{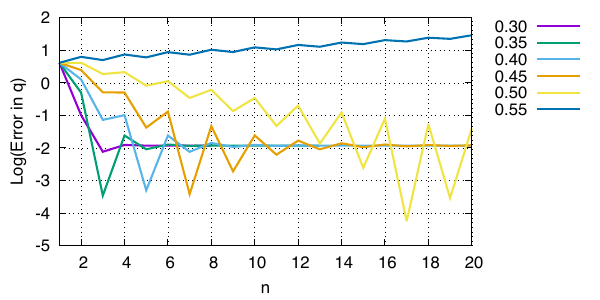} \\ 
\caption{Convergence history of the HDG solutions by the Neumann-Neumann algorithm.
The $L^2$-errors in $\qq$ are plotted for $\theta = 0.05, \ldots, 0.25$ (top) and 
$\theta = 0.30, \ldots, 0.55$ (bottom).} 
\label{fig-hdg-nn2-qerr}
\end{figure}

\begin{figure}[h]
\includegraphics[width=300pt]{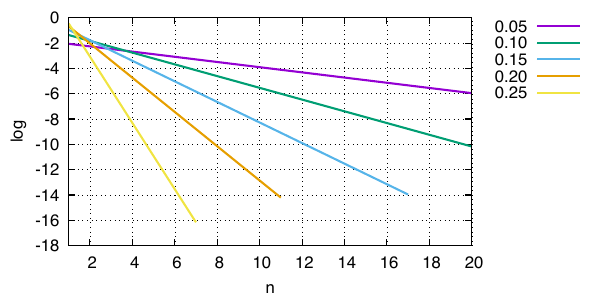}  \\
\includegraphics[width=300pt]{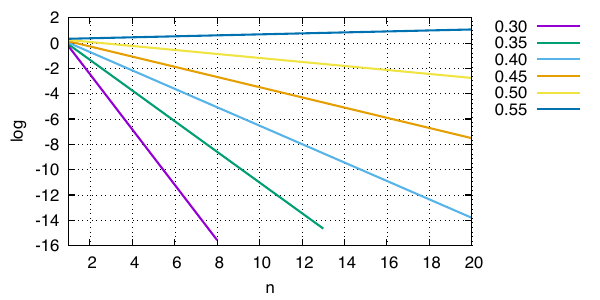}  
\caption{Difference $\|\widehat{u}_\Gamma^{n+1} - \widehat{u}_\Gamma^{n}\|_{L^2(\Gamma_{12})}$  are plotted in log scale for $\theta = 0.05, \ldots, 0.25$ (top) and 
$\theta = 0.30, \ldots, 0.55$ (bottom).} 
\label{fig-hdg-nn2-ugamma}
\end{figure}

The initial value is taken as $\widehat{u}_{\Gamma}^0 \equiv 0$.
The termination criteria is $\|\widehat{u}_{\Gamma}^{n+1} - \widehat{u}_{\Gamma}^{n}\|_{L^2(\Gamma_{12})} < 10^{-6}$.
We carried out numerical computation for $\theta = 0.05, 0.10, \ldots, 0.55$,
and histories of convergence in $\qq$ are displayed in Figures \ref{fig-hdg-nn2-qerr} and \ref{fig-hdg-nn2-ugamma}.
We do not show the errors of $u_i$ because they are very similar to the results of $\qq_i$.
We observe that our algorithm is convergent if $0 < \theta \le 0.5$ and 
the convergence speed is fastest around $\theta = 0.25$, 
which is similar to the case of the Neumann-Neumann method. 

\section{Trace--Flux alternating algorithm}
\subsection{Two-subdomain case} 
Let $u_\Gamma^n$ be a given function defined on the interface $\Gamma_{12}$.
We solve the subproblems with trace-interface condition
\begin{subequations} \label{traceif}
\begin{align}
-\Delta u_i^{n+1/2} &= f \text{ in } \Omega_i,  \\ 
u_i^{n+1/2} &= 0 \text{ on } \partial\Omega_i \cap \partial\Omega,  \\
u_i^{n+1/2} &= u_\Gamma^n \text{ on } \Gamma_{12}. 
\end{align} 
\end{subequations}
Then, we define an interface flux by
\[
\lambda_{\Gamma}^{n+1/2} = \frac{1}{2}\left(\pdv{u_1}{n_1}^{n+1/2}+\pdv{u_2}{n_2}^{n+1/2}\right).
\]
Solving the subproblems with flux-interface condition
\begin{subequations} \label{fluxif}
\begin{align}
-\Delta u_i^{n+1} &= f \text{ in } \Omega_i,  \\ 
u_i^{n+1} &= 0 \text{ on } \partial\Omega_i \cap \partial\Omega,  \\
\pdv{u_i}{n_i}^{n+1} &= (-1)^i\lambda_\Gamma^{n+1/2} \text{ on } \Gamma_{12},
\end{align}
\end{subequations}
we get $u_i^{n+1}$. The interface trace is updated by 
\[
u_\Gamma^{n+1} = \frac{1}{2}(u_1^{n+1} + u_2^{n+1}) \text{ on } \Gamma_{12}.
\]
Iteratively updating the interface trace by the above procedure, 
we can obtain the solution of the problem \eqref{pois} if $u_\Gamma^{n}$ converges.  

\subsection{HDG approximation}
The HDG approximation of the trace-flux alternating algorithm presented in the previous subsection 
is described as follows. 
\begin{enumerate}
\item[Step 1.] Set $\widehat{u}_\Gamma^0$.
\item[Step 2.] Repeat Steps 3-6 for $n \ge 0$ until convergence.
\item[Step 3.] Solve the trace-interface subproblems: Find $(\qq_{i}^{n+1/2}, u_{i}^{n+1/2}, \widehat{u}_{i}^{n+1/2}) \in \bm{V}_{ih}\times W_{ih} \times M_{ih}$ for $i=1,2$ such that
\begin{align*}
\intelem[ih]{\qq_i^{n+1/2}, \vv} - \intfa[ih]{u_i^{n+1/2}, \nabla\cdot \vv} 
+ \intfa[ih]{\widehat{u}^{n+1/2}, \vv\cdot\nn} & = 0 && \forall \vv \in V_{ih} \\
-\intelem[ih]{\qq_i^{n+1/2}, \nabla w} + \intfa[ih]{\widehat{\qq}^{n+1/2}\cdot\nn, w}  
&= (f,w)_{\Omega_i} && \forall w \in W_{ih} \\ 
\intfa[ih]{\widehat{\qq}^{n+1/2}\cdot\nn, \mu} & = 0 && \forall \mu \in M_h, 
\end{align*}
where 
\begin{align*}
\widehat{\qq}^{n+1/2}_i \cdot\nn &= \qq_i^{n+1/2} \cdot \nn + \tau(u_i^{n+1/2} - \widehat{u}^{n+1/2}_i), \\
\widehat{u}_i^{n+1/2} &= \widehat{u}_\Gamma^n  ~~\text{ on } \Gamma_{12}.
\end{align*}
\item[Step 4.] Define an interface flux by 
\[\lambda_{12}^{n+1/2} =
\frac{1}{2}(\widehat{\qq}_1^{n+1/2} \cdot \nn_1
+\widehat{\qq}_2^{n+1/2} \cdot \nn_2) \text{ on } \Gamma_{12}.
\]
\item[Step 5.] Solve the flux-interface subproblems: Find $(\qq_{i}^{n+1}, u_{i}^{n+1}, \widehat{u}_{i}^{n+1})\in \bm{V}_{ih}\times W_{ih} \times M_{ih}$ for $i=1,2$
such that
\begin{align*}
\intelem[ih]{\qq_i^{n+1}, \vv} - \intfa[ih]{u_i^{n+1}, \nabla\cdot \vv} 
+ \intfa[ih]{\widehat{u}^{n+1}, \vv\cdot\nn} & = 0 && \forall \vv \in V_{ih}, \\
-\intelem[ih]{\qq_i^{n+1}, \nabla w} + \intfa[ih]{\widehat{\qq}^{n+1}\cdot\nn, w}  
&= (f,w)_{\Omega_i} && \forall w \in W_{ih}, \\ 
\intfa[ih]{\widehat{\qq}^{n+1}\cdot\nn, \mu} & = (-1)^{i-1}\lambda_{12}^{n+1/2} && \forall \mu \in M_{ih}. 
\end{align*}
\item[Step 6.] Update the interface trace by 
\[ 
\widehat{u}_{12}^{n+1} = 
\frac{1}{2}(\widehat{u}_{1}^{n+1}
+\widehat{u}_{2}^{n+1}) \text{ on } \Gamma_{12}.
\]
\end{enumerate}

\subsection{For many-subdomain cases}  
\begin{figure}[h]
\centering
\includegraphics{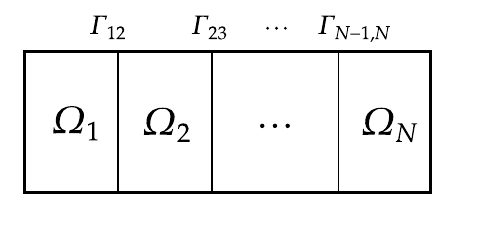}
\caption{Illustration of subdomains and interfaces}
\label{fig:manysubdomains}
\end{figure}

\begin{figure}[b]
\centering \includegraphics{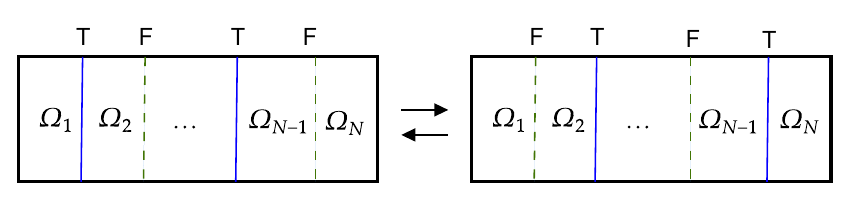}
\caption[width=200pt]{Illustration of the types of interfaces. Letters T and F means the types of trace- and 
flux-interfaces, respectively.}
\label{fig:interface-type}
\end{figure}

Let $\Omega$ be a disjoint union of $\Omega_1,\ldots, \Omega_{N}$ and 
define $\Gamma_{ij} = \partial\Omega_i \cap \partial\Omega_j$.
We here assume that $\overline{\Omega_i} \cap \overline{\Omega_j} =\emptyset$
if $|i - j| \ge 2$ and the length or area of $\Gamma_{i,i+1}$ is a positive value for $1 \le i \le N-1$, see Figure \ref{fig:manysubdomains}. 
We introduce types of interface to $\Gamma_{ij}$, which takes either trace- or flux-interface.
We propose the following algorithm.
\begin{enumerate}
\item[Step 1.] Set the types of $\Gamma_{12}, \Gamma_{34}, \ldots, \Gamma_{N-1,N}$ to be
trace-interfaces and set the others to be flux-interface.
\item[Step 2.] Set $\widehat{u}_{i,i+1}^0$ and $\lambda_{i,i+1}^0$ for $1 \le i \le N-1$.
\item[Step 3.] Repeat Steps 4-8 until convergence.
\item[Step 4.] Solve the following subproblems to get $u_i^{n+1/2}$ for $1 \le i \le N$:
\begin{align*}
-\Delta u_{i}^{n+1/2} &= f \text{ in } \Omega_{i}, && \\ 
u_{i}^{n+1/2} &= 0 \text{ on } \partial\Omega_i \cap \partial\Omega, && \\
u_{i}^{n+1/2} &=   \widehat{u}_{ij}^n \quad  \text{ on $\Gamma_{ij}$ if $\Gamma_{ij}$ is a trace-interface}, && \\ 
\pdv{u_{i}}{n_i}^{n+1/2} &= (-1)^i\lambda^n_{ij} \quad \text{ on $\Gamma_{ij}$ if $\Gamma_{ij}$ is a flux-interface}.&& 
\end{align*}
Here and in what follows, $j \in \{i-1, j+1\}$.
\item[Step 5.] Update the interface trace and flux by
\begin{align*}
\widehat{u}_{ij}^{n+1} &= \frac{1}{2}(\widehat{u}_i^{n+1/2} + \widehat{u}_j^{n+1/2})
\quad\text{ on  $\Gamma_{ij}$ if $\Gamma_{ij}$ is a flux-interface}, \\ 
\lambda_{ij}^{n+1} &= \frac{1}{2}\left(\pdv{u_{i}}{n_i}^{n+1/2} + \pdv{u_{j}}{n_j}^{n+1/2} \right)
\quad \text{ on  $\Gamma_{ij}$ if $\Gamma_{ij}$ is a trace-interface.}
\end{align*}
\item[Step 6.] Flip the types of interfaces. 
If the type of $\Gamma_{ij}$ is flux-interface, then set $\Gamma_{ij}$ to be trace-interface.
Else, if  the type of $\Gamma_{ij}$ is trace-interface, then set $\Gamma_{ij}$ to be flux-interface. 
See also Figure \ref{fig:interface-type}.
\item[Step 7.] Solve the following to get $u_i^{n+1}$:
\begin{align*}
-\Delta u_{i}^{n+1} &= f \text{ in } \Omega_{i}, && \\ 
u_{i}^{n+1} &= 0 \text{ on } \partial\Omega_i \cap \partial\Omega, && \\
u_{i}^{n+1} &=   \widehat{u}_{ij}^n \quad  \text{ on $\Gamma_{ij}$ if  the type of $\Gamma_{ij}$ is trace-interface}, && \\ 
\pdv{u_{i}}{n_i}^{n+1} &= (-1)^i\lambda^n_{ij} \quad \text{ on $\Gamma_{ij}$ if the type of $\Gamma_{ij}$ is flux-interface}. &&
\end{align*}
\item[Step 8.] Update the trace and flux 
\begin{align*}
\widehat{u}_{ij}^{n+1} &= \frac{1}{2}(\widehat{u}_i^{n+1/2} + \widehat{u}_j^{n+1/2})
\quad\text{ on  $\Gamma_{ij}$ if the type of $\Gamma_{ij}$ is  flux-interface}, \\ 
\lambda_{ij}^{n+1} &= \frac{1}{2}\left(\pdv{u_{i}}{n_i}^{n+1/2} + \pdv{u_{j}}{n_j}^{n+1/2} \right)
\quad \text{ on  $\Gamma_{ij}$ if the type of $\Gamma_{ij}$ is  trace-interface.}
\end{align*}
\end{enumerate}
This algorithm can be discretized by the HDG method in the same manner as in the two-subdomain case.
\subsection{Numerical results} 
In this section, we show the numerical results of 
the trace-flux algorithm for the test problem \eqref{testprob}.

\subsubsection{Two-subdomain case}
We study the dependence of convergence speed on the sizes of subdomains. 
Let $\alpha \in (0,0.5)$ and decompose $\Omega$ into
$\Omega_1 = (0,\alpha)\times(0,1)$ and $\Omega_2 = (\alpha,1)\times(0,1)$. 
We use unstructured meshes where $\Omega_1$ and $\Omega_2$
are divided into about $32\alpha \times 32$ and $32(1-\alpha)\times 32$ triangles, respectively,
and piecewise polynomials of degree 1.

\begin{figure}[h]
\includegraphics{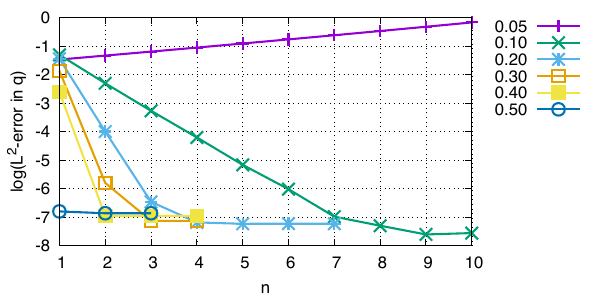}
\caption{$L^2$-errors in $\qq$ are plotted in log scale.}
\label{fig-tf2-aa}
\end{figure}

We computed solutions $(\qq_i, u_i, \widehat{u}_i)$ ~$(i=1,2)$ with $\alpha$ varying from $0.05$ to $0.5$
in order to study how the convergence property 
depends on the sizes of subdomains.
The history of convergence for various $\alpha$ is displayed in Figure \ref{fig-tf2-aa}.
When $\theta = 0.5$, i.e. the sizes of the subdomains are equal, the 
the iteration is terminated in 3 iterations and the convergence is fastest. 
As the parameter $\theta$ tends to zero, it takes more iterations to converge. 
For $\theta = 0.1, \ldots, 0.5$, the errors are monotonically decreasing and 
the final errors are similar.
In the case of $\theta = 0.05$, the errors are monotonically increasing and 
the solution seems to diverge. 
These results suggest that the convergence gets faster as subdomains get closer to uniform.

\subsubsection{Many-subdomain case}
\begin{figure}[h]
\centering 
\includegraphics{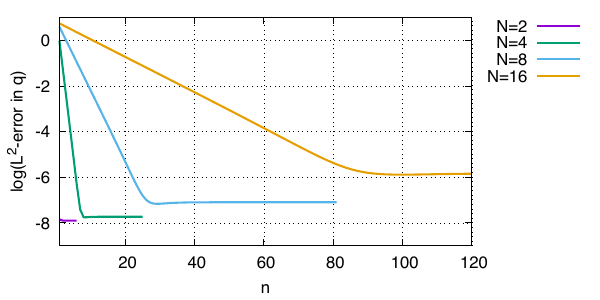}
\caption{$L^2$-errors in $\qq$  are plotted in log scale for various $N$}
\label{trace-flux-many}
\end{figure}

By numerical experiments, we demonstrate that our algorithm is valid for many-subdomain cases and 
examine its convergence property. 
The domain is equally divided into $N$ subdomains, and the width of 
a subdomain is $W = 1/N$. 
The $i$-th subdomain is denoted by $\Omega_i = ((i-1)W,iW)$ for $i \le i \le N$.
Unstructured meshes whose mesh size is $1/128$ and 
piecewise polynomials of degree 1 are used.
Figure \ref{trace-flux-many} shows the convergence history of the HDG solutions.
The solutions converge in 2, 16, 64, 240 iterations for 
$N = 2, 4, 8, 16$, respectively.
We see that the convergence gets slower as 
the number of division $N$ increase and its order is about $O(N^2)$.

\bibliographystyle{abbrv}
\bibliography{ref}
\end{document}